\documentclass{article}

\usepackage{microtype}
\usepackage{graphicx}
\usepackage{subfigure}
\usepackage{booktabs} 

\usepackage[accepted]{icml2019}

\icmltitlerunning{Simple Stochastic Gradient Methods for Non-Smooth Non-Convex Regularized 
Optimization}

\usepackage{amsmath,amssymb,amsthm}
\usepackage{tikz,pgfplots}

\def\ZZ{\mathbb{Z}}
\def\RR{\mathbb{R}}
\def\EE{\mathbb{E}}
\def\h{\tilde{h}}
\def\LL{{\cal L}}
\def\DD{\Delta}

\DeclareMathOperator*{\dist}{dist}

\newcommand{\prox}{\operatorname{prox}}

\newcommand{\argmin}{\operatornamewithlimits{argmin}}

\newtheorem{theorem}{Theorem}
\newtheorem{lemma}[theorem]{Lemma}
\newtheorem{property}[theorem]{Property}

\newtheorem{corollary}[theorem]{Corollary}

\begin{document}

\twocolumn[
\icmltitle{Simple Stochastic Gradient Methods for\\Non-Smooth Non-Convex Regularized Optimization}



\icmlsetsymbol{equal}{*}

\begin{icmlauthorlist}
\icmlauthor{Michael R. Metel}{ri}
\icmlauthor{Akiko Takeda}{ri,to}
\end{icmlauthorlist}

\icmlaffiliation{ri}{RIKEN Center for Advanced Intelligence Project, Tokyo, Japan}
\icmlaffiliation{to}{Department of Creative Informatics, Graduate School of Information Science 
and Technology, the University of Tokyo, Tokyo, Japan}

\icmlcorrespondingauthor{Michael R. Metel}{\mbox{michaelros.metel@riken.jp}}

\icmlkeywords{stochastic gradient, non-convex, non-smooth, regularization, optimization, variance 
	reduction}

\vskip 0.3in
]



\printAffiliationsAndNotice{}  

\begin{abstract}
Our work focuses on stochastic gradient methods for optimizing a smooth non-convex loss 
function with a non-smooth non-convex regularizer. Research on this class of problem is 
quite limited, and until recently no non-asymptotic convergence results have been reported. 
We present two simple stochastic gradient algorithms, for finite-sum and general stochastic 
optimization problems, which have superior convergence complexities compared to the current 
state-of-the-art. We also compare our algorithms' performance in practice for empirical risk 
minimization.    
\end{abstract}

\section{Introduction}
In this work we consider regularized optimization problems of the form
\begin{alignat}{6}\label{eq:1}
&\min\limits_{w\in\RR^d}&&\text{ }h(w):=f(w)+g(w),
\end{alignat}
where $f(w)$ has a Lipschitz continuous gradient and $g(w)$ has a proximal operator that can be 
efficiently computed. In addition, we assume that 
\begin{alignat}{6}\label{eq:genobj}
&f(w):=\EE_{\xi}[F(w,\xi)],
\end{alignat}
where $\xi\in \RR^p$ is a random vector following a probability distribution $P$ from which 
i.i.d. samples can be generated. We will also consider what is known as the finite-sum problem, 
where the expectation of $F(w,\xi)$ is taken over an empirical distribution function created by 
taking $n$ samples of $\xi$, $\xi_j$ for $j=1,...,n$: 
\begin{alignat}{6}\label{eq:fsobj}
&f(w):=\frac{1}{n}\sum_{j=1}^n f_j(w),
\end{alignat}
where $f_j(w)=F(w,\xi_j)$ and has a Lipschitz continuous gradient.

Our motivation for studying this problem is empirical risk minimization in machine learning. 
The purpose of $g(w)$, as a regularizer, is to induce a sparse solution when minimizing $f(w)$. 
Non-convex regularizers have been shown to outperform their convex counterparts with reduced bias 
in parameter estimation, including smoothly clipped absolute deviation (SCAD) \cite{fan2001} and 
minimax concave penalty (MCP) \cite{zhang2010}, as well as possess enhanced sparse signal 
recovery, such as the log-sum penalty \cite{candes2008}. In addition, improved generalization 
accuracy has been found using non-convex instead of convex loss functions \cite{shen2011}, with 
better robustness to outliers and noisy sample data \cite{wu2007, chapelle2009}. Smooth 
non-convex loss functions exhibiting these beneficial qualities include the sigmoid loss, Lorenz 
loss \cite{barbu22017}, and Savage loss \cite{masnadi2009}.

The literature concerning first-order stochastic methods for regularized optimization is vast, so 
we restrict our attention to algorithms achieving non-asymptotic rates of convergence for a 
non-convex function $f(w)$. Stochastic gradient methods for the case of a convex regularizer has 
been an active research area where algorithms with non-asymptotic convergence results were first 
achieved in \cite{ghadimi2016}. For finite-sum problems, \citet{reddi2016pro} 
were the first to develop a proximal algorithm using the stochastic variance reduced gradient 
approach of \citet{johnson2013}. The current state-of-the-art for the finite-sum problem seems to 
be the work of \citet{li2018} where one can also find a table of the convergence complexities of 
competing algorithms.

In the pursuit of solving \eqref{eq:1} where neither function $f(w)$ nor $g(w)$ are convex, 
the current body of research is quite limited. A generalization of \cite{ghadimi2016} with 
$g(w)$ being quasi-convex can be found in \cite{kawashima2018}, where the same convergence 
complexity is achieved. The only other work for non-convex regularizers to our knowledge is that 
of \citet{xu2018}, which recently improved upon the stochastic difference of convex (DC) algorithm 
of \citet{nitanda2017}, considering an objective of the form $c^1(w)-c^2(w)+g(w)$ where 
$c^1(w):=\EE_{\xi}[C^1(w,\xi)]$ and $c^2(w):=\EE_{\varsigma}[C^2(w,\varsigma)]$ are convex functions. 
It is assumed that $c^1(w)$ has a Lipschitz continuous gradient and $c^2(w)$ has a H\"older continuous 
gradient, and the proximal mapping of $g(w)$ can be efficiently computed. In their algorithms, a 
sequence of subproblems must be solved with increasing accuracy using a first-order stochastic 
algorithm, where convergence to a nearly $\epsilon$-critical point in a finite number of 
iterations is proved. The best convergence complexities in their work are achieved when it is 
assumed that $g(w)$ is Lipschitz continuous and $c^2(w)$ has a Lipschitz continuous gradient, 
which we will assume when discussing their work.\\

We now summarize the two main contributions of this paper:
\begin{itemize}
	\item Two algorithms are presented, a mini-batch stochastic gradient 
	algorithm for general stochastic objectives of the form \eqref{eq:genobj}, 
	and a variance reduced stochastic gradient algorithm for finite-sum 
	problems of the form \eqref{eq:fsobj}. We are aware of only 
	one other work, \cite{xu2018}, which has proven non-asymptotic convergence 
	for the class of problem we focus on in this paper. We attain superior 
	convergence results under both 
	objective assumptions, which are summarized in Table \ref{t:1}. The 
	complexities are 
	in terms of the number of gradient calls and proximal operations, see 
	Section 		
	\ref{prelims}.	
	\item No numerical experiments were conducted in \cite{xu2018}. We 
	implemented all algorithms for an application in empirical risk 
	minimization and found the simplest algorithm to 
	implement also performed the best in practice. 
\end{itemize} 

\begin{table*}[t]
	\caption{Comparison of convergence complexities obtained in \cite{xu2018} and this 
	paper.}
	\label{t:1}
	\begin{center}
			\begin{tabular}{llccc}
				\hline
				Algorithm	& Reference 
				&\begin{tabular}{@{}c@{}}Finite-sum\\Assumption\end{tabular}  & 
				\begin{tabular}{@{}c@{}}Gradient Call\\ 
					Complexity\end{tabular} & \begin{tabular}{@{}c@{}}Proximal Operator\\ 
					Complexity\end{tabular}\\
				\hline
				SSDC-SPG & \begin{tabular}{@{}c@{}}Theorem 7 a,\\\citet{xu2018}\end{tabular}
				& $\times$ & $O(\epsilon^{-8})$ & $O(\epsilon^{-8})$ \\
				SSDC-SVRG & \begin{tabular}{@{}c@{}}Theorem 7 c,\\\citet{xu2018}\end{tabular} & 
				$\surd$  & $O(n\epsilon^{-4})$ & $O(\epsilon^{-4})$\\
				MBSGA & Corollary \ref{mbcomp} & $\times$ & $O(\epsilon^{-5})$ &$O(\epsilon^{-4})$\\
				VRSGA & Corollary \ref{fscomp} & $\surd$ & $O(n^{2/3}\epsilon^{-3})$ 
				&$O(\epsilon^{-3})$\\
				\hline
			\end{tabular}
	\end{center}
	\vskip -0.1in
\end{table*}

\textbf{Remark:} In a subsequent revision uploaded after submission of this 
work, \citet{xu20182} 
present improved complexity results, as well as numerical experiments. The 
first row of Table 
\ref{t:1} would be $O(\epsilon^{-5})$ and $O(\epsilon^{-5})$, and the second 
row would be 
$\tilde{O}(n\epsilon^{-3})$ and 	$\tilde{O}(\epsilon^{-3})$ following the 
latest version of their 
work. 

\section{Preliminaries}
\label{prelims}
We assume that $f(w)$ has a Lipschitz continuous gradient with parameter $L$, 
$$||\nabla f(w)-\nabla f(x)||_2\leq L||w-x||_2,$$
which we will denote as being an $L$-smooth function. In the finite-sum case, we 
assume that each $f_j(w)$ is also $L$-smooth. Given a sample $\xi^k\sim P$, generated in iteration $k$ 
of an algorithm, we assume we can generate an unbiased stochastic gradient $\nabla F(w,\xi^k)$ such 
that  
\begin{alignat}{6}
&\EE[\nabla F(w,\xi^k)]=\nabla f(w),\label{eq:unbgrad}
\end{alignat}
and for some constant $\sigma$,
\begin{alignat}{6}
&\EE||\nabla F(w,\xi^k)-\nabla f(w)||^2_2\leq \sigma^2.\label{vbound}
\end{alignat}
Let $\partial h(w)$ denote the limiting subdifferential of our objective, defined as 
$$\partial h(w):=\{v: \exists w^k \xrightarrow{h} w, v^k\in \hat{\partial}h(w^k)\text{ with } 
v^k\rightarrow v\},$$
where $\hat{\partial}h(w):=\{v: \liminf\limits_{x\rightarrow w,x\neq w} \frac{h(x)-h(w)-\langle 
	v,x-w\rangle}{||x-w||_2}\geq 0\}$, and $w^k\xrightarrow{h}w$ signifies $w^k\rightarrow w$ with 
$h(w^k)\rightarrow h(w)$. The limiting subdifferential coincides with the gradient and subdifferential 
when the function is continuously differentiable and proper convex, respectively. We make use of the 
property that   
\begin{alignat}{6}
&\partial h(w)=\nabla f(w)+\partial g(w),\label{subprop}
\end{alignat}
for finite $g(w)$ \citep[Exercise 8.8 (c)]{rockafellar2009}. 
We also assume the proximal operator of $g(w)$ is nonempty for all $w\in \RR^d$ and $\lambda>0$, and 
can be efficiently computed, 
$$\prox_{\lambda g}(w):=\argmin\limits_{x\in 
	\RR^d}\left\{\frac{1}{2\lambda}||w-x||^2_2+g(x)\right\}.$$ 
In particular, let us denote an element as  
\begin{alignat}{6}
&\zeta^{\lambda}(w)\in \prox_{\lambda g}(w).\label{zprox}
\end{alignat}
We are interested in the convergence complexity of finding an $\epsilon$-stationary solution, such 
that for an algorithm solution $\overline{w}$,  
\begin{alignat}{6}
&\EE\left[\dist(0,\partial h(\overline{w}))\right]\leq \epsilon.\label{epsol}
\end{alignat}
We will measure algorithm complexity in terms of the number of gradient calls and proximal operations. 
For any $w$, a gradient call is  
either computing $\nabla F(w,\xi^k)$ given a sample $\xi^k$, or in the finite-sum case, returning 
$\nabla f_j(w)$ for a given $j$.

\section{Auxiliary functions of $h(w)$}
Our convergence results rely on bounding the gradient of a sequence of majorant functions 
of the auxiliary function  
$$\h_{\lambda}(w):=f(w) + e_{\lambda}g(w)$$
in expectation, where 
$$e_{\lambda}g(w):=\inf_{x\in \RR^d} \left\{\frac{1}{2\lambda}||w-x||^2_2+g(x)\right\}$$
is the Moreau envelope of $g(w)$. By considering $x=w$, we observe that 
\begin{alignat}{6}
&&&e_{\lambda}g(w)\leq g(w).\label{melb}
\end{alignat}
The Moreau envelope can be written as a DC function, 
\begin{alignat}{6}
&e_{\lambda}g(w)&&=\frac{1}{2\lambda}||w||^2_2-D^{\lambda}(w),\label{fsg}
\end{alignat}
where $D^{\lambda}(w)=\sup_{x\in 
	\RR^d}\left(\frac{1}{\lambda}w^Tx-\frac{1}{2\lambda}||x||^2_2-g(x)\right)$. 
We note that as the supremum of a set of affine functions, $D^{\lambda}(w)$ is convex, and 
we see from \eqref{zprox} that $\zeta^{\lambda}(w)$ attains the supremum of $D^{\lambda}(w)$. We can 
write down a smooth majorant of $\h_{\lambda}(w)$ as 
$$E_{\lambda}^k(w):=f(w)+U^k_{\lambda}(w)$$
in iteration $k$, where\\$U^k_{\lambda}(w)=\frac{1}{2\lambda}||w||^2_2
-\left(D^{\lambda}(w^k)+\frac{1}{\lambda}\zeta^{\lambda}(w^k)^T(w-w^k)\right)$.

The gradient of $E_{\lambda}^k(w)$ is 
\begin{alignat}{6}
&\nabla E_{\lambda}^k(w)=\nabla f(w)+\frac{1}{\lambda}(w-\zeta^{\lambda}(w^k)).\label{grade}
\end{alignat}
\begin{property}
	\label{Eprop}
	The following holds for $E^k_{\lambda}(w)$.		
	\begin{alignat}{6}
	&E_{\lambda}^k(w)\geq \h_{\lambda}(w)\hspace{0 mm} \text{ for all }w\in\RR^d\label{eq:e1}&\\
	&E_{\lambda}^k(w^k)=\h_{\lambda}(w^k)\label{eq:e3}&\\
	&E_{\lambda}^k(w) \text{ is } 
	L_{E\lambda}:=\left(L+\frac{1}{\lambda}\right)-\text{smooth.}&\label{eq:e2}
	\end{alignat}
\end{property}
\begin{proof}
	Given that both functions contain $f(w)$, it is sufficient to show that 
	\eqref{eq:e1} and \eqref{eq:e3} hold between the second terms $U^k_{\lambda}(w)$ and 
	$e_{\lambda}g(w)$.	
	
	\eqref{eq:e1}: As found in \cite{liu2017}, for any $w,z\in \RR^d$, 	
	\begin{alignat}{6}
	&&&D^{\lambda}(w)-D^{\lambda}(z)\nonumber\\
	&=&&\sup_{x\in 
		\RR^d}\left(\frac{1}{\lambda}w^Tx-\frac{1}{2\lambda}||x||^2-g(x)\right)\nonumber\\
	&&&-\sup_{x\in 
		\RR^d}\left(\frac{1}{\lambda}z^Tx-\frac{1}{2\lambda}||x||^2-g(x)\right)\nonumber\\
	&\geq&&\frac{1}{\lambda}w^T\zeta^{\lambda}(z)-\frac{1}{2\lambda}||\zeta^{\lambda}(z)||^2
	-g(\zeta^{\lambda}(z))\nonumber\\	
	&&&-\left(\frac{1}{\lambda}z^T\zeta^{\lambda}(z)-\frac{1}{2\lambda}||\zeta^{\lambda}(z)||^2-g(\zeta^{\lambda}(z))\right)\nonumber\\
	&=&&\frac{1}{\lambda}\zeta^{\lambda}(z)(w-z).\nonumber
	\end{alignat}	
	Setting $z=w^k$,  	
	\begin{alignat}{6}
	&e_{\lambda}g(w)&&=\frac{1}{2\lambda}||w||^2-D^{\lambda}(w)\nonumber\\
	&&&\leq\frac{1}{2\lambda}||w||^2-(D^{\lambda}(w^k)+\frac{1}{\lambda}\zeta^{\lambda}(w^k)^T(w-w^k))\nonumber\\
	&&&=U^k_{\lambda}(w).\nonumber
	\end{alignat}
			
	\eqref{eq:e3}:  
	$U^k_{\lambda}(w^k)=\frac{1}{2\lambda}||w^k||^2_2-D^{\lambda}(w^k)=e_{\lambda}g(w^k)$ from 
	\eqref{fsg}.
	
	\eqref{eq:e2}: \begin{alignat}{6}
	&&&\left\lVert\nabla E_{\lambda}^k(w)-\nabla E_{\lambda}^k(w')\right\lVert_2\nonumber\\
	&=&&\lVert\nabla f(w)+\frac{1}{\lambda}\left(w-\zeta^{\lambda}(w^k)\right)\nonumber\\
	&&&-\left(\nabla	
	f(w')+\frac{1}{\lambda}\left(w'-\zeta^{\lambda}(w^k)\right)\right)\lVert_2\nonumber\\
	&\leq&&(L+\frac{1}{\lambda})\lVert w-w'\rVert_2.\nonumber
	\end{alignat}	
\end{proof}
We note that the Moreau envelope of a convex function is also $\frac{1}{\lambda}$-smooth 
\citep[Theorem 6.60]{beck2017}, so there is no increase in the smoothness parameter for 
non-convex functions by taking a first-order approximation of the Moreau envelope.

\section{Mini-batch stochastic gradient algorithm}

\begin{algorithm}[t]
	\caption{Mini-batch stochastic gradient algorithm\\(MBSGA)}
	\label{alg:one}
	\begin{algorithmic}
		\STATE {\bfseries Input:} $w^{1}\in \RR^d$, $N\in\ZZ_{>0}$, $\alpha,\theta\in \RR$ 		
		\STATE $M:=\lceil N^{\alpha}\rceil\text{, }\lambda=\frac{1}{N^{\theta}}$
		\STATE $L_{E\lambda}=L+\frac{1}{\lambda}$
		\STATE $\gamma=\min\left\{\frac{1}{L_{E\lambda}},\frac{1}{\sigma\sqrt{N}}\right\}$
		\STATE $R\sim \text{uniform}\{1,...,N\}$
		\FOR{$k=1,2,...,R-1$} 
		\STATE $\zeta^{\lambda}(w^k)\in \prox_{\lambda g}(w^k)$
		\STATE Sample $\xi^k\sim P^M$
		\STATE $\nabla A^k_{\lambda M}(w^k,\xi^k)=\frac{1}{M}\sum_{j=1}^{M}\nabla 
		F(w^k,\xi^k_j)+\frac{1}{\lambda}(w^k-\zeta^{\lambda}(w^k))$
		\STATE $w^{k+1}=w^k-\gamma\nabla A^k_{\lambda M}(w^k,\xi^k)$
		\ENDFOR
		\STATE {\bfseries Output:} $\bar{w}^R\in \prox_{\lambda g}(w^R)$
	\end{algorithmic}
\end{algorithm}

\subsection{Convergence analysis}

The convergence analysis of MBSGA follows the technique of \citet{ghadlan2013} 
adapted to our problem. The following lemma bounds $\EE||\nabla E^R_{\lambda}(w^R)||^2_2$, 
with which we will ultimately bound $\EE\left[\dist(0,\partial h(\bar{w}^R)\right]$ in Theorem 
\ref{mbconv}.

\begin{lemma}	
	\label{th:1}
	For an initial value $w_1\in \RR^d$, $N\in\ZZ_{>0}$, and $\alpha,\theta \in \RR$, MBSGA 
	generates $w^R$ satisfying the following bound.  	
	\begin{alignat}{6} 
	&\EE||\nabla E^R_{\lambda}(w^R)||^2_2&&\leq 
	\frac{\tilde{\DD}}{N}(L+N^{\theta})	
	+\frac{\sigma}{\sqrt{N}}\left(\tilde{\DD}+\frac{L+N^{\theta}}{\lceil 
		N^{\alpha}\rceil}\right),\nonumber
	\end{alignat}	
	where $\tilde{\DD}=2(\h_{\lambda}(w^1)-\h_{\lambda}(w^*_{\lambda}))$ and 
	$w^*_{\lambda}$ is a global minimizer of $\h_{\lambda}(\cdot)$.
\end{lemma}

Due to a lack of space, the proof of Lemma \ref{th:1} can be found in Section 1 of the 
supplementary material. In order to prove the convergence of $\EE\left[\dist(0,\partial 
h(\bar{w}^R)\right]$, we will 
require the following two properties.
\begin{property}	
	\label{eq:switch}
	Assume that $g(w)$ is Lipschitz continuous with parameter $l$,
	\begin{alignat}{6}
	&&\dist(0,\partial h(\zeta^{\lambda}(w^k)))&\leq ||\nabla 
	E^k_{\lambda}(w^k)||_2+2l\lambda 
	L\nonumber.
	\end{alignat}	
\end{property}
\begin{proof}
	Given that $\zeta^{\lambda}(w)$ is a minimizer of\\$\frac{1}{2\lambda}||w-x||^2_2+g(x)$ from 
	\eqref{zprox}, $$\frac{1}{\lambda}(w-\zeta^{\lambda}(w))\in \partial g(\zeta^{\lambda}(w))$$ 
	and $$\nabla f(\zeta^{\lambda}(w^k))+\frac{1}{\lambda}(w^k-\zeta^{\lambda}(w^k))\in \partial 
	h(\zeta^{\lambda}(w^k))$$	
	using \eqref{subprop}. 
	It follows that
	\begin{alignat}{6}
	&&&\dist(0,\partial h(\zeta^{\lambda}(w^k)))\nonumber\\
	&\leq&& ||\nabla 
	f(\zeta^{\lambda}(w^k))+\frac{1}{\lambda}(w^k-\zeta^{\lambda}(w^k))||_2\nonumber\\
	&=&&||\nabla f(w^k)-\nabla f(w^k)+\nabla 
	f(\zeta^{\lambda}(w^k))\nonumber\\
	&&&+\frac{1}{\lambda}(w^k-\zeta^{\lambda}(w^k))||_2\nonumber\\
	&\leq&&||\nabla f(w^k)+\frac{1}{\lambda}(w^k-\zeta^{\lambda}(w^k))||_2\nonumber\\
	&&&+||\nabla f(\zeta^{\lambda}(w^k))-\nabla f(w^k)||_2\nonumber\\
	&\leq&&||\nabla 
	E^k_{\lambda}(w^k)||_2+L||w^k-\zeta^{\lambda}(w^k)||_2.\nonumber
	\end{alignat}		
	In order to bound $||w^k-\zeta^{\lambda}(w^k)||_2$, recall from 
	\eqref{melb} that    				  
	\begin{alignat}{6}
	&&g(w)&\geq e_{\lambda}g(w)\nonumber\\
	&&&=\frac{1}{2\lambda}||w-\zeta^{\lambda}(w)||_2^2+g(\zeta^{\lambda}(w)).\nonumber
	\end{alignat}
	Rearranging and using the Lipschitz continuity assumption,
	\begin{alignat}{6}
	&&\frac{1}{2\lambda}||w-\zeta^{\lambda}(w)||_2^2&\leq 
	g(w)-g(\zeta^{\lambda}(w))\nonumber\\
	&&&\leq l||w-\zeta^{\lambda}(w)||_2\nonumber\\
	&&||w-\zeta^{\lambda}(w)||_2&\leq 2l\lambda.\nonumber 
	\end{alignat}	
\end{proof}
\begin{property}	
	\label{le:2}	
	Let $w^*$ be a global minimizer of $h(\cdot)$ and let $w^*_{\lambda}$ be a global minimizer 
	of $\h_{\lambda}(\cdot)$. Assume that $g(w)$ is Lipschitz continuous with 
	parameter $l$, then 	
	$$\h_{\lambda}(w)-\h_{\lambda}(w^*_{\lambda})\leq h(w)-h(w^*) + 
	\frac{l^2\lambda}{2}.$$	
\end{property}
\begin{proof}	
	\begin{alignat}{6}
	&&&\h_{\lambda}(w)-\h_{\lambda}(w^*_{\lambda})-\left(h(w)-h(w^*)\right)\nonumber\\
	&=&&e_{\lambda}g(w) - f(w_{\lambda}^*) -e_{\lambda}g(w_{\lambda}^*)\nonumber\\
	&&&-\left(g(w) - f(w^*)-g(w^*)\right)\nonumber\\
	&\leq&&- f(w_{\lambda}^*) -e_{\lambda}g(w_{\lambda}^*)+f(w^*) +g(w^*)\nonumber\\
	&\leq&& - f(w_{\lambda}^*) -e_{\lambda}g(w_{\lambda}^*)+f(w_{\lambda}^*) 
	+g(w_{\lambda}^*)\nonumber\\
	&=&&g(w_{\lambda}^*)-e_{\lambda}g(w_{\lambda}^*),\nonumber
	\end{alignat}	
	where the first inequality follows from \eqref{melb}. For any $w$, by the definition of the 
	Moreau envelope,
	\begin{alignat}{6}
	&&e_{\lambda}g(w)&=\frac{1}{2\lambda}||w-\zeta^{\lambda}(w)||_2^2+g(\zeta^{\lambda}(w))\nonumber\\
	&&g(w)-e_{\lambda}g(w)&=g(w)-g(\zeta^{\lambda}(w))-\frac{1}{2\lambda}||w-\zeta^{\lambda}(w)||_2^2\nonumber\\
	&&&\leq l||w-\zeta^{\lambda}(w)||_2-\frac{1}{2\lambda}||w-\zeta^{\lambda}(w)||_2^2.\nonumber
	\end{alignat}	
	The right-hand side is maximized when $||w-\zeta^{\lambda}(w)||_2=l\lambda$, giving the 
	desired result, 		
	\begin{alignat}{6}\label{eq:propK1} &g(w)-e_{\lambda}g(w)\leq \frac{l^2\lambda}{2}.
	\end{alignat}	
\end{proof}
We note that \eqref{eq:propK1} cannot be 
improved under the further assumption that $g(w)$ is convex, which can be found in \citep[Theorem 
10.51]{beck2017}.
\begin{theorem}	
	\label{mbconv}	
	Assume that $g(w)$ is Lipschitz continuous with parameter $l$. The output $\bar{w}^R$ of 
	MBSGA satisfies  	
	\begin{alignat}{6}
	&\EE\left[\dist(0,\partial h(\bar{w}^R))\right]\leq 
	\sqrt{\frac{(\DD+l^2N^{-\theta})(L+N^{\theta})}{N}}\nonumber\\
	&+\sqrt{\frac{\sigma}{\sqrt{N}}\left(\DD+\frac{l^2}{N^{\theta}}+\frac{L+N^{\theta}}{\lceil 
			N^{\alpha}\rceil}\right)}+\frac{2lL}{N^{\theta}},\nonumber
	\end{alignat}
	where $\DD=2(h(w^1)-h(w^*))$ and $w^*$ is a global minimizer of $h(\cdot)$.	
\end{theorem}
\begin{proof}
	From Property \ref{eq:switch}, choosing $\zeta^{\lambda}(w^R)=\bar{w}^R$,
	$${\dist(0,\partial h(\bar{w}^R))\leq ||\nabla 	
	E^R_{\lambda}(w^R)||_2+2l\lambda L}.$$ 
	Taking its expectation, 	
	\begin{alignat}{6}
	&&&\EE\left[\dist(0,\partial h(\bar{w}^R))\right]\nonumber\\
	&\leq&&\EE[||\nabla E^R_{\lambda}(w^R)||_2]+2l\lambda	L\nonumber\\
	&\leq&&\sqrt{\EE\left[||\nabla 
	E^R_{\lambda}(w^R)||^2_2\right]}+\frac{2lL}{N^{\theta}} 
	\nonumber\\
	&\leq&&\sqrt{\frac{\tilde{\DD}(L+N^{\theta})}{N}}+\sqrt{\frac{\sigma}{\sqrt{N}}\left(\tilde{\DD}+\frac{L+N^{\theta}}{\lceil
			N^{\alpha}\rceil}\right)}+\frac{2lL}{N^{\theta}},\nonumber
	\end{alignat}	
	where the second inequality follows from Jensen's inequality and the third inequality uses Lemma 
	\ref{th:1}.
	The result then follows using Property \ref{le:2} as  	
	\begin{alignat}{6}
	&&\tilde{\DD}=2(\h_{\lambda}(w^1)-\h_{\lambda}(w^*_{\lambda}))&\leq 2(h(w^1)-h(w^*)) + 
	l^2\lambda\nonumber\\
	&&&=\DD+\frac{l^2}{N^{\theta}}\nonumber
	\end{alignat} 
	
\end{proof}

Now that we have bounded the expected distance of $\partial h(\bar{w}^R)$ from the origin, we 
prove an $\epsilon$-stationary point convergence complexity.

\begin{corollary}	
	\label{mbcomp}
	Assume that $g(w)$ is Lipschitz continuous with 
	parameter $l$. To obtain an $\epsilon$-stationary solution \eqref{epsol} using MBSGA, the 
	gradient call complexity is $O(\epsilon^{-5})$ and the proximal operator complexity is 
	$O(\epsilon^{-4})$ when $\alpha=\theta=0.25$. 	
\end{corollary}
\begin{proof}	
	From Theorem \ref{mbconv},	
	\begin{alignat}{6}
	&&&\EE\left[\dist(0,\partial h(\bar{w}^R))\right]\nonumber\\
	&\leq&&\sqrt{\frac{(\DD+l^2N^{-\theta})(L+N^{\theta})}{N}}\nonumber\\
	&&&+\sqrt{\frac{\sigma}{\sqrt{N}}\left(\DD+\frac{l^2}{N^{\theta}}+\frac{L+N^{\theta}}{\lceil 
	N^{\alpha}\rceil}\right)}+\frac{2lL}{N^{\theta}}\nonumber\\
	&&=&O(N^{0.5\theta-0.5})+O(N^{-0.25}+N^{0.5\theta-0.5\alpha-0.25})\nonumber\\
	&&&+O(N^{-\theta}).\nonumber
	\end{alignat}
	
	Setting $\theta=\alpha=0.25$,  	
	\begin{alignat}{6}
	&&\EE\left[\dist(0,\partial h(\bar{w}^R))\right]&\leq 
	O(N^{-0.25}).\nonumber
	\end{alignat}
	
	An $\epsilon$-stationary solution will require less than $N=O(\epsilon^{-4})$ 
	iterations. One proximal operation is done per iteration, which establishes the proximal operator 
	complexity of $O(\epsilon^{-4})$. The number of gradient calls per iteration is $\lceil 
	N^\alpha\rceil=O(\epsilon^{-1})$. The number of gradient calls to get an 
	$\epsilon$-stationary solution is then less than	
	\begin{alignat}{6}
	&N\lceil N^\alpha\rceil=O(\epsilon^{-5}).\nonumber
	\end{alignat}	
\end{proof}
\section{Variance reduced method for finite-sum problems}
In this section we assume that 
$$f(w)=\frac{1}{n}\sum_{j=1}^nf_j(w),$$
where each $f_j(w)$ is $L$-smooth. 
\begin{algorithm}[t]
	\caption{Variance reduced stochastic gradient algorithm (VRSGA)}
	\label{alg:two}
	\begin{algorithmic}
		\STATE {\bfseries Input:} $\tilde{w}^1\in \RR^d$, $N\in\ZZ_{>0}$, $\alpha,\theta\in \RR$
		\STATE $m=\lceil n^{\alpha}\rceil$, $b=m^2$ 
		\STATE $S=\lceil\frac{N}{m}\rceil$, $\lambda=(Sm)^{-\theta}$
		\STATE $L_{E\lambda}=L+\frac{1}{\lambda}$, $\gamma=\frac{1}{6L_{E\lambda}}$
		\STATE $R\sim \text{uniform}\{1,...,S\}$
		\FOR{$k=1,2,...,R$}
		\STATE $w^{k}_1=\tilde{w}^{k}$
		\STATE $G^k=\nabla f(\tilde{w}^{k})$
		\FOR{$t=1,2,...,m$}
		\STATE $\zeta^{\lambda}(w^k_t)\in \prox_{\lambda g}(w^k_t)$
		\STATE $I\sim \text{uniform}\{1,...,n\}^b$
		\STATE $V^k_t=\frac{1}{b}\sum_{j\in I}\left(\nabla f_j(w^k_t)-\nabla 		
		f_j(\tilde{w}^k)\right)+G^k+\frac{1}{\lambda}(w^k_t-\zeta^{\lambda}(w^k_t))$
		\STATE $w^k_{t+1}=w^k_t-\gamma V^k_t$
		\ENDFOR		
		\STATE $\tilde{w}^{k+1}=w^k_{m+1}$
		\ENDFOR
		\STATE $T\sim \text{uniform}\{1,...,m\}$		
		\STATE {\bfseries Output:} $\bar{w}^R_T\in \prox_{\lambda g}(w^R_T)$
	\end{algorithmic}
\end{algorithm}

\subsection{Convergence analysis}
In our convergence analysis, we make use of the function $E^k_{t\lambda}(w)$, which is 
constructed in the same manner as $E^k_{\lambda}(w)$, using $w^k_t$ instead of $w^k$. This 
function possesses the same characteristics as found in Property \ref{Eprop}. The convergence 
analysis follows closely to the work of \citet{li2018} adapted to our problem.
\begin{lemma}	
	\label{vrlemma}
	For an initial value $\tilde{w}_1\in \RR^d$, $N\in\ZZ_{>0}$, and $\alpha,\theta\in \RR$, 
	VRSGA generates $w^R_T$ satisfying the following bound.	  	  	
	\begin{alignat}{6} 
	&&\EE\left[||\nabla 
	E^R_{T\lambda}(w^R_T)||^2_2\right]&\leq\tilde{\DD}\frac{L+(Sm)^{\theta}}{Sm},\nonumber
	\end{alignat}	
	where $\tilde{\DD}=36(\h_{\lambda}(\tilde{w}^1)-\h_{\lambda}(w^*_{\lambda}))$ and 
	$w^*_{\lambda}$ is a global minimizer of $\h_{\lambda}(\cdot)$.	
\end{lemma}
The proof of Lemma \ref{vrlemma} can be found in Section 2 of the supplementary material. 
\begin{theorem}	
	\label{vrconv}	
	Assume that $g(w)$ is Lipschitz continuous with parameter $l$. The output $\bar{w}^R_T$ of 
	VRSGA satisfies 
	\begin{alignat}{6}
	&&&\EE\left[||\dist(0,\partial 
	h(\bar{w}^R_T))||_2\right]\nonumber\\
	&\leq&&\sqrt{\frac{\left(L+(Sm)^{\theta}\right) 
			\left(\DD+18l^2(Sm)^{-\theta}\right)}{Sm}}+\frac{2lL}{(Sm)^{\theta}},\nonumber
	\end{alignat}		
	where $\DD=36(h(w^1)-h(w^*))$ and $w^*$ is a global minimizer of $h(\cdot)$.		
\end{theorem}
\begin{proof}
	The proof follows what was done to prove Theorem \ref{mbconv}. From Property 
	\ref{eq:switch}, 
	$${\dist(0,\partial h(\bar{w}^R_T))\leq ||\nabla 
	E^R_{T\lambda}(w^R_T)||_2+2l\lambda L}.$$ 
	Taking its expectation, 	
	\begin{alignat}{6}
	&&&\EE\left[||\dist(0,\partial h(\bar{w}^R_T))||_2\right]\nonumber\\
	&\leq&&\EE[||\nabla E^R_{T\lambda}(w^R_T)||_2]+2l\lambda L\nonumber\\
	&\leq&&\sqrt{\EE\left[||\nabla 
	E^R_{T\lambda}(w^R_T)||^2_2\right]}+\frac{2lL}{(Sm)^{\theta}} 
	\nonumber\\
	&\leq&&\sqrt{\frac{\left(L+(Sm)^{\theta}\right) 
			\tilde{\DD}}{Sm}}+\frac{2lL}{(Sm)^{\theta}}\nonumber\\
	&\leq&&\sqrt{\frac{\left(L+(Sm)^{\theta}\right)			
			\left(\DD+18l^2(Sm)^{\theta}\right)}{Sm}}+\frac{2lL}{(Sm)^{\theta}},\nonumber
	\end{alignat}
	where the third inequality follows from Lemma \ref{vrlemma}. The fourth inequality holds 
	using Property \ref{le:2},   	
	\begin{alignat}{6}
	&&\tilde{\DD}=36(\h_{\lambda}(w^1)-\h_{\lambda}(w^*_{\lambda}))&\leq 36(h(w^1)-h(w^*)) + 
	18l^2\lambda\nonumber\\
	&&&=\DD+\frac{18l^2}{(Sm)^{\theta}}.\nonumber
	\end{alignat} 		
\end{proof}
\begin{corollary}	
	\label{fscomp}
	Assume that $g(w)$ is Lipschitz continuous with parameter $l$. To obtain an $\epsilon$-stationary 
	solution \eqref{epsol} using VRSGA, the gradient call complexity is 
	$O(n^{\frac{2}{3}}\epsilon^{-3})$ and the proximal operator complexity is $O(\epsilon^{-3})$ 
	choosing $\alpha=\theta=\frac{1}{3}$.
\end{corollary}
\begin{proof}
	From Theorem \ref{vrconv} with $\theta=\frac{1}{3}$,	
	\begin{alignat}{6}
	&&&\EE\left[||\dist(0,\partial h(\bar{w}^R_T))||_2\right]\nonumber\\
	&\leq&&\sqrt{\frac{\left(L+(Sm)^{\frac{1}{3}}\right) 
			\left(\DD+18l^2(Sm)^{\frac{-1}{3}}\right)}{Sm}}+\frac{2lL}{(Sm)^{\frac{1}{3}}}\nonumber\\
	&&&=O((Sm)^{-\frac{1}{3}})\nonumber
	\end{alignat}	
	
	An $\epsilon$-stationary solution will require at most  $Sm=O(\epsilon^{-3})$ 
	iterations, which establishes the proximal operator complexity. The number of gradient calls after 
	$Sm$ iterations, taking $\alpha=\frac{1}{3}$ is $$Sn+Smb=Sm\frac{n}{\lceil 
		n^{\frac{1}{3}}\rceil}+Sm\lceil n^{\frac{1}{3}}\rceil^2=O(n^{\frac{2}{3}}\epsilon^{-3}).$$	
\end{proof}

\section{Application}
\label{app}
In this section we consider the application of binary classification for a particular choice of loss 
function and regularizer, which will be used in our numerical experiments. Non-convex Lipschitz 
continuous regularizers which have proximal operators with closed form solutions include the log-sum 
penalty, SCAD, MCP, and the capped $l_1$-norm. For their closed form solutions, see \cite{gong2013}. 
All of these functions are separable, 
$g(w):=\sum_{i=1}^dg_i(w_i)$. For $\kappa, \nu>0$, the log-sum penalty is
$$g_i(w_i)=\kappa\log(1+|w_i|/\nu).$$
\begin{property}
	\label{firstprop}	
	The log-sum penalty is $\frac{\kappa}{\nu}\sqrt{d}$-Lipschitz continuous.
\end{property}
\begin{proof}	
	Assume $w_i\geq 0$ over which $g_i(w_i)$ is differentiable and 
	$|\frac{dg_i}{dw_i}(w_i)|\leq 
	\frac{\kappa}{\nu}$. Using the mean value theorem with $z_i\geq 0$, 
	$|g_i(z_i)-g_i(w_i)|\leq 
	\frac{\kappa}{\nu}|z_i-w_i|$. Given the symmetry of $g_i(w_i)$, this bound 
	holds for general $w_i$ 
	and $z_i$. It then follows that for any $w$ and $z$,
	\begin{alignat}{6}
	&&|g(z)-g(w)|&=\left|\sum_{i=1}^d(g_i(z_i)-g_i(w_i))\right|\nonumber\\
	&&&\leq \sum_{i=1}^d|g_i(z_i)-g_i(w_i)|\nonumber\\
	&&&\leq \frac{\kappa}{\nu}\sum_{i=1}^d|z_i-w_i|\nonumber\\
	&&&\leq \frac{\kappa}{\nu}\sqrt{d}||z-w||_2\nonumber
	\end{alignat}		
\end{proof}
Smooth non-convex loss functions, which are known to be robust to outliers, include the sigmoid 
loss, $\frac{1}{1+e^v}$, Lorenz loss 
\citep{barbu22017}, Savage loss \citep{masnadi2009}, and the tangent loss \citep{masnadi2010}. 
We will consider the Lorenz loss,
$$\LL(v)=\begin{cases}
0 & \text{ if } v>1 \\ 
\log(1+(v-1)^2)  & \text{ otherwise }\\ 
\end{cases}$$ 
for $v\in\RR$, which is differentiable everywhere \citep{barbu22017}. For the problem setting of 
binary classification, we have a set of training data 
$\{x,y\}$ where 
$y=\{y^1,y^2,...,y^n\}$, $y^j\in\{-1,1\}$, is the label set, and $x=\{x^1,x^2,...,x^n\}$, $x^j\in 
\RR^d$, is the feature set. Our loss function is then
$$f(w)=\frac{1}{n}\sum_{j=1}^n f_j(w),$$
where
$$f_j(w)=\LL(y^jw^Tx^j).$$
\begin{property}	
	\label{secprop}
	Using the Lorenz loss function, $f(w)$ is $\frac{2}{n}\sum_{j=1}^n||x^j||^2_2$-smooth.	
\end{property}
\begin{proof}	
	We first consider the function 
	$$\hat{\LL}(v)=\log(1+(v-1)^2).$$	
	Its first and second derivatives are
	$$\hat{\LL}'(v)=\frac{2(v-1)}{1+(v-1)^2}$$
	and
	$$\hat{\LL}''(v)=\frac{2}{1+(v-1)^2}-\left(\frac{2(v-1)}{1+(v-1)^2}\right)^2.$$
		
	We can see that $v=1$ maximizes $\hat{\LL}''(v)$, with $\hat{\LL}''(1)=2$. 
	Examining the third 
	derivative, 
	$$\hat{\LL}'''(v)=\frac{-4(v-1)}{\left(1+(v-1)^2\right)^2}\left(3-\frac{4(v-1)^2}{1+(v-1)^2}\right),$$
	$v=1\pm \sqrt{3}$ minimizes $\hat{\LL}''(v)$ with $\hat{\LL}''(1\pm 
	\sqrt{3})=-0.25$, so we 
	conclude that 
	$$|\hat{\LL}''(v)|\leq|\hat{\LL}''(1)|=2.$$	
	Using the mean value theorem, for any $v$ and $u$, 	
	\begin{alignat}{6}
	&&|\hat{\LL}'(v)-\hat{\LL}'(u)|&\leq 2|v-u|.\nonumber
	\end{alignat}	
	We now show that $\LL(v)$ is also 2-smooth. For $v>1$, 
	$\LL'(v)=\hat{\LL}'(1)=0$. Taking $v>1$ and 
	$u\leq 1$, 
	\begin{alignat}{6}
	&|\LL'(v)-\LL'(u)|&&=|\hat{\LL}'(1)-\hat{\LL}'(u)|\nonumber\\
	&&&\leq 2 |1-u|\nonumber\\
	&&&\leq 2 |v-u|.\nonumber
	\end{alignat}	
	An $L$-smooth function composed with the linear function, $y^jw^Tx^j$, is 
	$L||y^jx^j||^2_2$-smooth \citep[Claim 12.9]{shalev2014}, so $f_j(w)$ is 
	$2||x^j ||^2_2$-smooth and the result follows.	
\end{proof}
We also note that the Lorenz loss function is DC-decomposable, which is required to implement the 
algorithms of \cite{xu2018}. 
\begin{property}	
	\label{lastprop}
	The Lorenz loss function is DC-decomposable, 	
	$$\LL(v)=\LL^1(v)-\LL^2(v),$$
	where $\LL^1(v)=\frac{1}{8}v^2+\LL(v)$ and $\LL^2(v)=\frac{1}{8}v^2$.
\end{property}
\begin{proof}
	Since $\LL''(v)\geq-\frac{1}{4}$, from the proof of Property \ref{secprop}, 
	we write the DC 
	decomposition of $\LL(v)$ as $\LL(v)=\LL^1(v)-\LL^2(v)$, where 
	$\LL^1(v)=\frac{1}{8}v^2+\LL(v)$ 
	and $\LL^2(v)=\frac{1}{8}v^2$.
\end{proof}

\section{Numerical experiments}
\label{NE}
We conducted experiments comparing our algorithms to those of \citep{xu2018} for the problem of binary 
classification as described in Section \ref{app}, on datasets a9a \citep{dua2017} and MNIST 
\citep{lecun1998}, as used in \citep{reddi2016pro,allen2016,li2018}. For the MNIST dataset, our 
objective was to learn class 1. The dimensions of a9a are $n=32,561$ and $d=123$, and those of MNIST 
are $n=60,000$ and $d=784$. All experiments were conducted using MATLAB 2017b on a Mac Pro with a 2.7 
GHz 12-core Intel Xeon E5 processor and 64GB of RAM. 
We compare performance in terms of the log of the objective function and 
wall-clock 
time. 

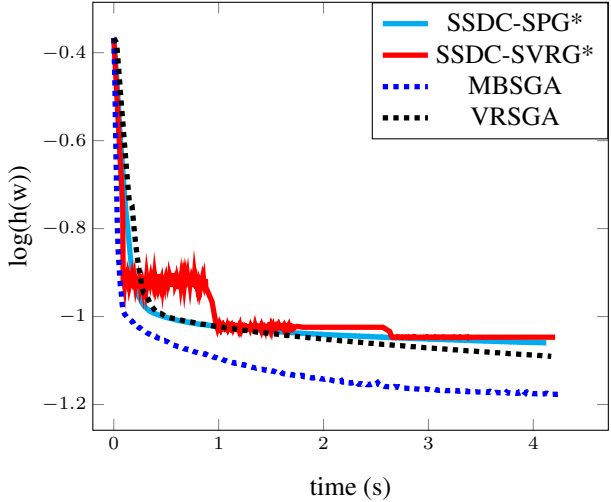
\begin{figure}[t]
	\begin{center}
		\centerline{
			\begin{tikzpicture}
			\scriptsize
			\begin{axis}[legend style={font=\normalsize},label 
			style={font=\normalsize},xlabel=time (s),ylabel=log(h(w)),xmin=-0.2,
			legend style={at={(1,1)}, anchor=north east},
			y label style={at={(axis description cs:0.075,.5)},anchor=south}]
			\addplot[line width=2pt,draw=cyan]
			table[x=x,y=ly]
			{XuSPGa.dat};
			\addplot[line width=2pt,draw=red]
			table[x=x,y=ly]
			{XuSVRGa.dat};	
			\addplot[line width=2pt,draw=blue,dotted]
			table[x=x,y=ly]
			{MBSGAa.dat};
			\addplot[line width=2pt,draw=black,dotted]
			table[x=x,y=ly]
			{SVRGAa.dat};
			\legend{SSDC-SPG*,SSDC-SVRG*,MBSGA,VRSGA}
			\end{axis}
			\end{tikzpicture}
		}
		\caption{Comparison of algorithms of this paper and \citep{xu2018} (marked with *) using 
			the a9a dataset} \label{T3}
	\end{center}
\end{figure}

All algorithms' convergence rates rely on outputting a random iteration. In 
order to fairly compare algorithms we ignore this step, e.g. for MBSGA, we set 
$R=N$. The algorithms were initially run taking $e=15$ effective passes over 
the data for a9a and $e=9$ for MNIST. These values were adjusted so that all 
algorithms ended at approximately the same time. The regularizer's 
parameters were chosen as $\kappa=\frac{1}{d}$ and $\nu=1$. All parameter 
values used in MBSGA and VRSGA were obtained from the theoretical convergence 
results, except for the upper bound $\sigma$ \eqref{vbound} used in MBSGA. This 
parameter was estimated by doing 50 
iterations of MBSGA with step size $\gamma=\frac{1}{L_{E\lambda}}$, using a 
different random seed than was used for the experiments, and computing the sample estimate 
$\hat{\sigma}^k$ each iteration with the $M$ samples used in the algorithm. An 
estimate of $\sigma$ was then taken as $\hat{\sigma}=\max_{k}\hat{\sigma}^k$. 

The proof of convergence of algorithms VRSGA and SSDC-SVRG rely on the 
assumption that each $f_j(w)$ is $L$-smooth, so for these instances 
$L=2\max_j||x^j||^2_2$. For algorithms MBSGA and VRSGA, the final proximal 
operation at the output was omitted and can be considered as simply a means of 
proving the non-asymptotic convergence of the algorithms. 

No experiments were done in \citep{xu2018}, so we implemented their algorithms following the parameter 
values found in their theoretical results and remarks, and recommended in 
\cite{xiao2014}, from which their work is partially based on. Full details of their algorithms' 
implementation can be found in Section 3 of the supplementary material.

Figures \ref{T3} and \ref{T4} show the results of the experiments. We observe 
that MBSGA outperformed all other algorithms. MBSGA is also the 
simplest algorithm to implement, making it an appealing choice for use in practice. It 
appears all other algorithms would require further parameter tuning in order for them to 
possibly perform comparably. 

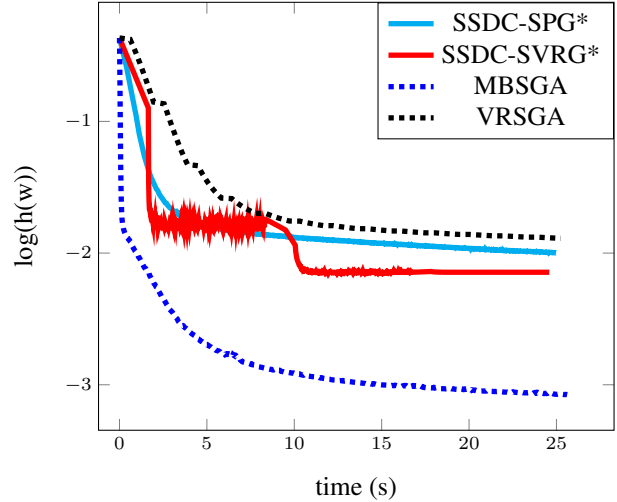
\begin{figure}[t]
	\begin{center}
		\centerline{
			\begin{tikzpicture}
			\scriptsize
			\begin{axis}[legend style={font=\normalsize},label 
			style={font=\normalsize},xlabel=time (s),ylabel=log(h(w)),
			legend style={at={(1,1)}, anchor=north east},xmin=-1.2,
			y label style={at={(axis description cs:0.075,.5)},anchor=south}]				
			\addplot[line width=2pt,draw=cyan]
			table[x=x,y=ly]
			{XuSPGm.dat};
			\addplot[line width=2pt,draw=red]
			table[x=x,y=ly]
			{XuSVRGm.dat};
			\addplot[line width=2pt,draw=blue,dotted]
			table[x=x,y=ly]
			{MBSGAm.dat};
			\addplot[line width=2pt,draw=black,dotted]
			table[x=x,y=ly]
			{SVRGAm.dat};			
			\legend{SSDC-SPG*,SSDC-SVRG*,MBSGA,VRSGA}
			\end{axis}
			\end{tikzpicture}
		}
		\caption{Comparison of algorithms of this paper and \citep{xu2018} (marked with *) using 
			the	MNIST dataset} 
		\label{T4}
	\end{center}
\end{figure}

\section{Conclusion and future research}
We have presented two simple stochastic gradient algorithms for optimizing a 
smooth non-convex loss function with a non-smooth non-convex regularizer. 
Our work improves upon the only other known non-asymptotic convergence results 
of \citet{xu2018} for this class of problem. Superior convergence complexities 
were shown for the case of a general stochastic loss 
function using a mini-batch stochastic gradient algorithm, and for the case of 
a finite-sum loss function using a variance reduced stochastic gradient algorithm. 
In an empirical study we found that the simplest algorithm to implement was 
also the best performing, making it the most appealing algorithm considered for 
this problem setting. Future research using the techniques developed in this 
work could consider additional regularizers in the objective to induce 
desirable properties of the solution in addition to sparsity.  

\section*{Acknowledgements}
This work was supported by JSPS KAKENHI Grant Numbers 15K00031 and 19H04069, and supported by
JST CREST Grant Numbers JPMJCR15K5 and JPMJCR14D2.

\raggedright
\bibliography{SSGM}
\bibliographystyle{icml2019}

\end{document}



\onecolumn

\icmltitle{Simple Stochastic Gradient Methods for Non-Smooth\\Non-Convex Regularized Optimization: 
	Supplementary Material}

\icmlsetsymbol{equal}{*}

\section{Proof of Lemma 2}
\begin{customlemma}{2}
	For an initial value $w_1\in \RR^d$, $N\in\ZZ_{>0}$, and $\alpha,\theta \in \RR$, MBSGA 
	generates $w^R$ satisfying the following bound.  	
	\begin{alignat}{6} 
	&\EE||\nabla E^R_{\lambda}(w^R)||^2_2&&\leq 
	\frac{\tilde{\DD}}{N}(L+N^{\theta})	
	+\frac{\sigma}{\sqrt{N}}\left(\tilde{\DD}+\frac{L+N^{\theta}}{\lceil 
		N^{\alpha}\rceil}\right),\nonumber
	\end{alignat}	
	where $\tilde{\DD}=2(\h_{\lambda}(w^1)-\h_{\lambda}(w^*_{\lambda}))$ and 
	$w^*_{\lambda}$ is a global minimizer of $\h_{\lambda}(\cdot)$.
\end{customlemma}
In order to prove this result, we require the following property.
\begin{property}	
	\label{eq:minib}	
	\begin{alignat}{6}
	&&\EE||\nabla A^k_{\lambda M}(w^k,\xi^k)-\nabla 
	E^k_{\lambda}(w^k)||^2_2&\leq\frac{\sigma^2}{M}\nonumber
	\end{alignat}	
\end{property}
\begin{proof}
	From the definition of $\nabla A^k_{\lambda M}(w^k,\xi^k)$ found in Algorithm 1 and 
	(11), $\nabla A^k_{\lambda M}(w^k,\xi^k)-\nabla 	
	E^k_{\lambda}(w^k)=\frac{1}{M}\sum_{j=1}^{M}\nabla F(w^k,\xi^k_j)-\nabla 
	f(w^k)$. Taking the expectation of its squared norm,
	\begin{alignat}{6}
	&&\EE||\nabla A^k_{\lambda M}(w^k,\xi^k)-\nabla 
	E^k_{\lambda}(w^k)||^2_2&=\EE||\frac{1}{M}\sum_{j=1}^{M}(\nabla F(w^k,\xi^k_j)-\nabla 
	f(w^k))||^2_2\nonumber\\
	&&&=\frac{1}{M^2}\EE\sum_{i=1}^n\left(\sum_{j=1}^M\nabla F(w^k,\xi^k_j)_i-\nabla 
	f(w^k)_i\right)^2.\nonumber
	\end{alignat}	
	For $j\neq l$, $\nabla F(w^k,\xi^k_j)_i-\nabla f(w^k)_i$ and $\nabla F(w^k,\xi^k_l)_i-\nabla 
	f(w^k)_i$ 
	are independent random variables with zero mean. It follows that    
	\begin{alignat}{6}
	&\EE[(\nabla F(w^k,\xi^k_j)_i-\nabla f(w^k)_i)(\nabla F(w^k,\xi^k_l)_i-\nabla 
	f(w^k)_i)]=\nonumber\\
	&\EE[(\nabla F(w^k,\xi^k_j)_i-\nabla f(w^k)_i)]\EE[(\nabla F(w^k,\xi^k_l)_i-\nabla 
	f(w^k)_i)]=0,&\nonumber
	\end{alignat}	
	and 	
	\begin{alignat}{6}
	&&\frac{1}{M^2}\EE\sum_{i=1}^n\left(\sum_{j=1}^M\nabla F(w^k,\xi^k_j)_i-\nabla 
	f(w^k)_i\right)^2&=\frac{1}{M^2}\EE\sum_{i=1}^n\sum_{j=1}^{M}(\nabla F(w^k,\xi^k_j)_i-\nabla 
	f(w^k)_i)^2\nonumber\\
	&&&=\frac{1}{M^2}\sum_{j=1}^{M}\EE||\nabla F(w^k,\xi^k_j)-\nabla 
	f(w^k)||^2_2\leq\frac{\sigma^2}{M}\nonumber
	\end{alignat}	
	using (5).
\end{proof}
\begin{proof}[Proof of Lemma 2]	
	Given the smoothness of $E^k_{\lambda}(w)$ as shown in Property 1,
	\begin{alignat}{6}
	&E^k_{\lambda}(w^{k+1})&&\leq E^k_{\lambda}(w^k)+\langle\nabla 
	E^k_{\lambda}(w^k),w^{k+1}-w^k\rangle+\frac{L_{E\lambda}}{2}||w^{k+1}-w^k||^2_2\nonumber\\ 
	&&&=E^k_{\lambda}(w^k)+\langle\nabla 
	E^k_{\lambda}(w^k),-\gamma\nabla A^k_{\lambda 
		M}(w^k,\xi^k)\rangle+\frac{L_{E\lambda}}{2}||-\gamma\nabla A^k_{\lambda 
		M}(w^k,\xi^k)||^2_2.\nonumber
	\end{alignat}	
	Using (12) and (13),	
	\begin{alignat}{6}
	&\h(w^{k+1})&&\leq \h(w^k)-\gamma\langle\nabla 
	E^k_{\lambda}(w^k),\nabla A^k_{\lambda 
		M}(w^k,\xi^k)\rangle+\frac{L_{E\lambda}}{2}\gamma^2||\nabla A^k_{\lambda 
		M}(w^k,\xi^k)||^2_2.\nonumber
	\end{alignat}	
	Setting $\delta_k=\nabla A^k_{\lambda M}(w^k,\xi^k)-\nabla E^k_{\lambda}(w^k)$, 
	\begin{alignat}{6}
	&&\h(w^{k+1})\leq& \h(w^k)-\gamma\left(||\nabla 
	E^k_{\lambda}(w^k)||^2_2+\langle\nabla 
	E^k_{\lambda}(w^k),\delta_k\rangle\right)+\frac{L_{E\lambda}}{2}\gamma^2\left(||\nabla 
	E^k_{\lambda}(w^k)||^2_2+2\langle\nabla 
	E^k_{\lambda}(w^k),\delta_k\rangle+||\delta_k||^2_2\right)\nonumber\\
	&&=&\h(w^k)+\left(\frac{L_{E\lambda}}{2} \gamma^2-\gamma\right)||\nabla 
	E^k_{\lambda}(w^k)||^2_2+(L_{E\lambda}\gamma^2-\gamma)\langle\nabla 
	E^k_{\lambda}(w^k),\delta_k\rangle+\frac{L_{E\lambda}}{2}\gamma^2||\delta_k||^2_2,\nonumber
	\end{alignat}	
	as		
	$$\langle\nabla E^k_{\lambda}(w^k),\nabla A^k_{\lambda M}(w^k,\xi^k)\rangle=||\nabla 
	E^k_{\lambda}(w^k)||^2_2+\langle\nabla E^k_{\lambda}(w^k),\delta_k\rangle$$
	and	
	$$||\nabla A^k_{\lambda M}(w^k,\xi^k)||^2_2=||\nabla E^k_{\lambda}(w^k)||^2_2+2\langle\nabla 
	E^k_{\lambda}(w^k),\delta_k\rangle+||\delta_k||^2_2.$$	
	After $N$ iterations, 	
	\begin{alignat}{6}
	\left(\gamma-\frac{L_{E\lambda}}{2}\gamma^2\right)\sum_{k=1}^N ||\nabla 
	E^k_{\lambda}(w^k)||^2_2&\leq&& \h(w^1)-\h(w^{N+1})+(L_{E\lambda} 
	\gamma^2-\gamma)\sum_{k=1}^N\langle\nabla 
	E^k_{\lambda}(w^k),\delta_k\rangle+\frac{L_{E\lambda}}{2}\gamma^2\sum_{k=1}^N||\delta_k||^2_2\nonumber\\
	&\leq&& \h_{\lambda}(w^1)-\h_{\lambda}(w^{*}_{\lambda})+(L_{E\lambda} 
	\gamma^2-\gamma)\sum_{k=1}^N\langle\nabla 
	E^k_{\lambda}(w^k),\delta_k\rangle+\frac{L_{E\lambda}}{2}\gamma^2\sum_{k=1}^N||\delta_k||^2_2.\nonumber
	\end{alignat}	
	It follows from (4) that for $w$ independent of $\xi^k$, $\EE \nabla 
	A^k_{\lambda M}(w,\xi^k)=\nabla E^k_{\lambda}(w)$, and so $\EE[\delta_k]=0$. Taking the 
	expectation 
	of both sides,	
	\begin{alignat}{6}
	&& \left(\gamma-\frac{L_{E\lambda}}{2}\gamma^2\right)\sum_{k=1}^N\EE||\nabla 	
	E^k_{\lambda}(w^k)||^2_2\leq& 
	\h(w^1)-\h(w^*_{\lambda})+\frac{L_{E\lambda}}{2}\gamma^2\sum_{k=1}^N\EE||\delta_k||^2_2\nonumber\\
	&&\leq& 
	\h(w^1)-\h(w^*_{\lambda})+\frac{L_{E\lambda}}{2}\gamma^2\frac{N}{M}\sigma^2,\nonumber
	\end{alignat}	
	where the second inequality uses Property \ref{eq:minib}. Choosing $R$ uniformly over 
	$\{1,...,N\}$,	
	\begin{alignat}{6}
	&&\EE||\nabla 
	E^R_{\lambda}(w^R)||^2_2&=\frac{1}{N}\sum_{k=1}^N\EE||\nabla E^k_{\lambda}(w^k)||^2_2\nonumber\\
	&&&\leq\frac{1}{N\left(\gamma-\frac{L_{E\lambda}}{2}\gamma^2\right)}\left(\h(w^1)-\h(w^*)
	+\frac{L_{E\lambda}}{2}\gamma^2\frac{N}{M}\sigma^2\right).\nonumber
	\end{alignat}	
	Since $\gamma\leq \frac{1}{L_{E\lambda}}$, it holds that 
	$\gamma-\frac{L_{E\lambda}}{2}\gamma^2\geq \frac{1}{2}\gamma$, and	
	\begin{alignat}{6}
	&&\frac{1}{N\left(\gamma-\frac{L_{E\lambda}}{2}\gamma^2\right)}\left(\frac{\tilde{\DD}}{2}+
	\frac{L_{E\lambda}}{2}\gamma^2\frac{N}{M}\sigma^2\right)\leq&
	\frac{1}{N\gamma}\left(\tilde{\DD}+L_{E\lambda}\gamma^2\frac{N}{M}\sigma^2\right)\nonumber\\
	&&=&\frac{\tilde{\DD}}{N\gamma} +L_{E\lambda} \frac{\gamma}{M}\sigma^2\nonumber\\
	&&\leq&\frac{\tilde{\DD}}{N}\max\left\{L_{E\lambda},\sigma\sqrt{N}\right\}+L_{E\lambda}\frac{\sigma}{M\sqrt{N}}\nonumber\\
	&&\leq&\frac{\tilde{\DD}L_{E\lambda}}{N}+\frac{\sigma}{\sqrt{N}}\left(\tilde{\DD}+\frac{L_{E\lambda}}{M}\right)\nonumber\\
	&&=&\frac{\tilde{\DD}}{N}(L+N^{\theta})	
	+\frac{\sigma}{\sqrt{N}}\left(\tilde{\DD}+\frac{L+N^{\theta}}{\lceil 
		N^{\alpha}\rceil}\right)\nonumber
	\end{alignat}
	
\end{proof}

\section{Proof of Lemma 7}
\begin{customlemma}{7}	
	For an initial value $\tilde{w}_1\in \RR^d$, $N\in\ZZ_{>0}$, $\alpha,\theta\in \RR$, VRSGA 
	generates $w^R_T$ satisfying the following bound.	  	  	
	\begin{alignat}{6} 
	&&\EE\left[||\nabla 
	E^R_{T\lambda}(w^R_T)||^2_2\right]&\leq\tilde{\DD}\frac{L+(Sm)^{\theta}}{Sm},\nonumber
	\end{alignat}	
	where $\tilde{\DD}=36(\h_{\lambda}(\tilde{w}^1)-\h_{\lambda}(w^*_{\lambda}))$ and 
	$w^*_{\lambda}$ is a global minimizer of $\h_{\lambda}(\cdot)$.	
\end{customlemma}
In order to prove this result, we require the following lemmas.
\begin{lemma}	
	\label{vrprop1}
	Consider arbitrary $w,V,z\in \RR^d$, $\gamma\in \RR$, and $w^+=w-\gamma V$,  	
	\begin{alignat}{6} 
	&E^k_{t\lambda}(w^+)&\leq& E^k_{t\lambda}(z)+\langle\nabla 
	E^k_{t\lambda}(w)-V,w^+-z\rangle+\frac{L_{E\lambda}}{2}||w^+-w||^2_2+\frac{L_{E\lambda}}{2}||z-w||^2_2-\frac{1}{\gamma}\langle
	 w^+-w,w^+-z\rangle.\nonumber	
	\end{alignat}
\end{lemma}
\begin{proof}
	Adding the following three inequalities proves the result, where the first two  
	come from the smoothness of $E^k_{t\lambda}(w)$ and $-E^k_{t\lambda}(w)$, see Property 	
	1, and the third is due to $V+\frac{1}{\gamma}(w^+-w)=0$.
	\begin{alignat}{6} &&E^k_{t\lambda}(w^+)&\leq E^k_{t\lambda}(w)+\langle\nabla 
	E^k_{t\lambda}(w),w^+-w\rangle+\frac{L_{E\lambda}}{2}||w^+-w||^2_2\nonumber\\
	&&-E^k_{t\lambda}(z)&\leq -E^k_{t\lambda}(w)+\langle-\nabla 
	E^k_{t\lambda}(w),z-w\rangle+\frac{L_{E\lambda}}{2}||z-w||^2_2\nonumber\\
	&&0&=-\langle V+\frac{1}{\gamma}(w^+-w),w^+-z\rangle\nonumber
	\end{alignat}
\end{proof}
\begin{lemma}	
	\label{yineq}
	For vectors $w$, $x$, $z$, and $\beta>0$, 	
	\begin{alignat}{6} 
	&&||w-x||^2_2&\leq 
	(1+\beta)||w-z||^2_2+\left(1+\frac{1}{\beta}\right)||z-x||^2_2.\nonumber
	\end{alignat} 	
	
\end{lemma}

\begin{proof}
	
	\begin{alignat}{6} 
	&&||w-x||^2_2&=||w-z+z-x||^2_2\nonumber\\
	&&&\leq \left(||w-z||_2+||z-x||_2\right)^2\nonumber\\
	&&&=||w-z||^2_2+2||w-z||_2||z-x||_2+||z-x||^2_2\nonumber\\
	&&&\leq||w-z||^2_2+\left(\beta||w-z||^2_2+\frac{1}{\beta}||z-x||^2_2\right)+||z-x||^2_2\nonumber\\
	&&&=(1+\beta)||w-z||^2_2+\left(1+\frac{1}{\beta}\right)||z-x||^2_2,\nonumber
	\end{alignat} 
	
	where the second inequality uses Young's inequality.	
	
\end{proof}

\begin{proof}[Proof of Lemma 7]
	
	Let $\hat{w}^k_{t+1}=w^k_t-\gamma \nabla E^k_{t\lambda}(w^k_t)$, with $w^+=w^k_{t+1}$, 
	$w=w^k_t$, $V=V^k_t$, and $z=\hat{w}^k_{t+1}$ in Lemma \ref{vrprop1} to get the inequality
	\begin{alignat}{6} 
	&E^k_{t\lambda}(w^k_{t+1})&\leq& E^k_{t\lambda}(\hat{w}^k_{t+1})+\langle\nabla 	
	E^k_{t\lambda}(w^k_t)-V^k_t,w^k_{t+1}-\hat{w}^k_{t+1}\rangle+\frac{L_{E\lambda}}{2}||w^k_{t+1}-w^k_t||^2_2\nonumber\\
	&&&+\frac{L_{E\lambda}}{2}||\hat{w}^k_{t+1}-w^k_t||^2_2-\frac{1}{\gamma}\langle 
	w^k_{t+1}-w^k_t,w^k_{t+1}-\hat{w}^k_{t+1}\rangle.\label{vr:ineq2}
	\end{alignat}
	
	In addition, let $w^+=\hat{w}^k_{t+1}$, $w=w^k_t$, $V=\nabla E^k_{t\lambda}(w^k_t)$, and 
	$z=w^k_t$ in Lemma \ref{vrprop1} to get	
	\begin{alignat}{6} 
	&E^k_{t\lambda}(\hat{w}^k_{t+1})&\leq& E^k_{t\lambda}(w^k_t)+\langle\nabla 
	E^k_{t\lambda}(w^k_t)-\nabla 
	E^k_{t\lambda}(w^k_t),\hat{w}^k_{t+1}-w^k_{t+1}\rangle+\frac{L_{E\lambda}}{2}||\hat{w}^k_{t+1}-w^k_t||^2_2\nonumber\\
	&&&+\frac{L_{E\lambda}}{2}||w^k_t-w^k_t||^2_2-\frac{1}{\gamma}\langle 
	\hat{w}^k_{t+1}-w^k_t,\hat{w}^k_{t+1}-w^k_t\rangle\nonumber\\	
	&&=&E^k_{t\lambda}(w^k_t)+\left(\frac{L_{E\lambda}}{2}-\frac{1}{\gamma}\right)||\hat{w}^k_{t+1}-w^k_t||^2_2.\label{vr:ineq3}
	\end{alignat}
	
	Adding \eqref{vr:ineq2} and \eqref{vr:ineq3}, 	
	\begin{alignat}{6} 
	&E^k_{t\lambda}(w^k_{t+1})&\leq& E^k_{t\lambda}(w^k_t)+\langle\nabla 	
	E^k_{t\lambda}(w^k_t)-V^k_t,w^k_{t+1}-\hat{w}^k_{t+1}\rangle+\frac{L_{E\lambda}}{2}||w^k_{t+1}-w^k_t||^2_2\nonumber\\
	&&&-\frac{1}{\gamma}\langle 
	w^k_{t+1}-w^k_t,w^k_{t+1}-\hat{w}^k_{t+1}\rangle+	
	\left(L_{E\lambda}-\frac{1}{\gamma}\right)||\hat{w}^k_{t+1}-w^k_t||^2_2.\label{keyineq}
	\end{alignat}
	
	Plugging $\langle w^k_{t+1}-w^k_t,w^k_{t+1}-\hat{w}^k_{t+1}\rangle=\frac{1}{2}\left(|| 
	w^k_{t+1}-w^k_t||^2_2+||w^k_{t+1}-\hat{w}^k_{t+1}||^2_2-||\hat{w}^k_{t+1}-w^k_t||^2_2\right)$ 
	into \eqref{keyineq} and rearranging,	
	\begin{alignat}{6} 
	&E^k_{t\lambda}(w^k_{t+1})&\leq& E^k_{t\lambda}(w^k_t)+\langle\nabla 	
	E^k_{t\lambda}(w^k_t)-V^k_t,w^k_{t+1}-\hat{w}^k_{t+1}\rangle+\left(\frac{L_{E\lambda}}{2}-\frac{1}{2\gamma}\right)||w^k_{t+1}-w^k_t||^2_2\nonumber\\
	&&&-\frac{1}{2\gamma}||w^k_{t+1}-\hat{w}^k_{t+1}||^2_2+	
	\left(L_{E\lambda}-\frac{1}{2\gamma}\right)||\hat{w}^k_{t+1}-w^k_t||^2_2.\label{byi}
	\end{alignat} 
	
	Focusing on the term $-\frac{1}{2\gamma}||w^k_{t+1}-\hat{w}^k_{t+1}||^2_2$, we apply Lemma 
	\ref{yineq} with $w=w^k_{t+1}$, $x=w^k_t$, and $z=\hat{w}^k_{t+1}$. Rearranging, 	
	\begin{alignat}{6} 
	&&-(1+\beta)||w^k_{t+1}-\hat{w}^k_{t+1}||^2_2&\leq 
	-||w^k_{t+1}-w^k_t||^2_2+\left(1+\frac{1}{\beta}\right)||\hat{w}^k_{t+1}-w^k_t||^2_2\nonumber\\
	&&-\frac{1}{2\gamma}||w^k_{t+1}-\hat{w}^k_{t+1}||^2_2&\leq 
	-\frac{1}{(1+\beta)2\gamma} 	
	||w^k_{t+1}-w^k_t||^2_2+\frac{\left(1+\frac{1}{\beta}\right)}{(1+\beta)2\gamma}||\hat{w}^k_{t+1}-w^k_t||^2_2.\nonumber
	\end{alignat} 
	Choosing $\beta=3$,
	\begin{alignat}{6} 
	&&-\frac{1}{2\gamma}||w^k_{t+1}-\hat{w}^k_{t+1}||^2_2&\leq 
	-\frac{1}{8\gamma} 
	||w^k_{t+1}-w^k_t||^2_2+\frac{1}{6\gamma}||\hat{w}^k_{t+1}-w^k_t||^2_2.\nonumber
	\end{alignat} 
	
	Using this inequality in \eqref{byi},
	\begin{alignat}{6} 
	&E^k_{t\lambda}(w^k_{t+1})&\leq& E^k_{t\lambda}(w^k_t)+\langle\nabla 
	E^k_{t\lambda}(w^k_t)-V^k_t,w^k_{t+1}-\hat{w}^k_{t+1}\rangle+\left(\frac{L_{E\lambda}}{2}-\frac{1}{2\gamma}\right)||w^k_{t+1}-w^k_t||^2_2\nonumber\\
	&&&-\frac{1}{8\gamma} ||w^k_{t+1}-w^k_t||^2_2+\frac{1}{6\gamma}||\hat{w}^k_{t+1}-w^k_t||^2_2+	
	\left(L_{E\lambda}-\frac{1}{2\gamma}\right)||\hat{w}^k_{t+1}-w^k_t||^2_2\nonumber\\
	&&=& E^k_{t\lambda}(w^k_t)+\langle\nabla 
	E^k_{t\lambda}(w^k_t)-V^k_t,w^k_{t+1}-\hat{w}^k_{t+1}\rangle+\left(\frac{L_{E\lambda}}{2}-\frac{5}{8\gamma}\right)||w^k_{t+1}-w^k_t||^2_2\nonumber\\
	&&&+\left(L_{E\lambda}-\frac{1}{3\gamma}\right)||\hat{w}^k_{t+1}-w^k_t||^2_2\nonumber\\
	&&=& E^k_{t\lambda}(w^k_t)+\gamma ||\nabla 
	E^k_{t\lambda}(w^k_t)-V^k_t||^2_2+\left(\frac{L_{E\lambda}}{2}-\frac{5}{8\gamma}\right)||w^k_{t+1}-w^k_t||^2_2
	+\left(L_{E\lambda}-\frac{1}{3\gamma}\right)||\hat{w}^k_{t+1}-w^k_t||^2_2,\nonumber
	\end{alignat}
	where the last equality holds since $w^k_{t+1}-\hat{w}^k_{t+1}=\gamma(\nabla 
	E^k_{t\lambda}(w^k_t)-V^k_t)$. Using (12) and (13), and taking the 
	expectation of both sides,	
	\begin{alignat}{6} 
	&\EE\h_{\lambda}(w^k_{t+1})&\leq& 
	\EE\left[\h_{\lambda}(w^k_t)+\gamma                    
	||\nabla E^k_{t\lambda}(w^k_t)-V^k_t||^2_2
	+\left(\frac{L_{E\lambda}}{2}-\frac{5}{8\gamma}\right)||w^k_{t+1}-w^k_t||^2_2
	+\left(L_{E\lambda}-\frac{1}{3\gamma}\right)||\hat{w}^k_{t+1}-w^k_t||^2_2\right].\label{bvb}
	\end{alignat}	
	Focusing on $\EE\left[||\nabla E^k_{t\lambda}(w^k_t)-V^k_t||^2_2\right]$,  	
	from (11) and the definition of $V^k_t$ found in Algorithm 2, $\nabla 
	E^k_{t\lambda}(w^k_t)-V^k_t=\nabla f(w^k_t)-(\frac{1}{b}\sum_{j\in 
		I}\left(\nabla f_j(w^k_t)-\nabla f_j(\tilde{w}^k)\right)+G^k)$. Rearranging, and taking the 
	expectation of its squared norm,	
	\begin{alignat}{6} 
	&&\EE||\nabla E^k_{t\lambda}(w^k_t)-V^k_t||^2_2&=\EE||\frac{1}{b}\sum_{j\in I}\left(\nabla 
	f_j(\tilde{w}^k)-\nabla f_j(w^k_t)\right)-\left(G^k-\nabla f(w^k_t)\right)||^2_2\nonumber\\
	&&&=\frac{1}{b^2}\EE\sum_{j\in I}||\nabla f_j(\tilde{w}^k)-\nabla f_j(w^k_t)-\left(G^k-\nabla 
	f(w^k_t)\right)||^2_2\nonumber\\
	&&&\leq\frac{1}{b^2}\EE\sum_{j\in I}||\nabla f_j(\tilde{w}^k)-\nabla 
	f_j(w^k_t)||^2_2\nonumber\\
	&&&\leq\frac{L^2}{b}\EE||\tilde{w}^k-w^k_t||^2_2.\nonumber
	\end{alignat}	
	As the squared norm of a sum of independent random variables with zero mean, the second 
	equality holds using the same reasoning as found in Property \ref{eq:minib}, and the first 
	inequality holds since $\EE||x-\EE[x]||^2_2\leq \EE||x||^2_2$ for any random variable $x$. Using 
	this bound in \eqref{bvb},	
	\begin{alignat}{6} 
	&\EE\h_{\lambda}(w^k_{t+1})&\leq& \EE\left[\h_{\lambda}(w^k_t)+\gamma 
	\frac{L^2}{b}||\tilde{w}^k-w^k_t||^2_2+
	\left(\frac{L_{E\lambda}}{2}-\frac{5}{8\gamma}\right)||w^k_{t+1}-w^k_t||^2_2+\left(L_{E\lambda}-\frac{1}{3\gamma}\right)||\hat{w}^k_{t+1}-w^k_t||^2_2\right]\nonumber\\
	&&\leq&\EE\left[\h_{\lambda}(w^k_t)+\frac{L_{E\lambda}}{6b}||\tilde{w}^k-w^k_t||^2_2-
	\frac{13L_{E\lambda}}{4}||w^k_{t+1}-w^k_t||^2_2-L_{E\lambda}||\hat{w}^k_{t+1}-w^k_t||^2_2\right]\nonumber\\
	&&=&\EE\left[\h_{\lambda}(w^k_t)+\frac{L_{E\lambda}}{6b}||\tilde{w}^k-w^k_t||^2_2-
	\frac{13L_{E\lambda}}{4}||w^k_{t+1}-w^k_t||^2_2-\frac{1}{36L_{E\lambda}}||\nabla	
	E^k_{t\lambda}(w^k_t)||^2_2\right],\label{byi2}
	\end{alignat}		
	where the last two lines use the fact that $\gamma=\frac{1}{6L_{E\lambda}}$. Focusing on 
	$-\frac{13L_{E\lambda}}{4}||w^k_{t+1}-w^k_t||^2_2$, we apply Lemma \ref{yineq} with 
	$w=w^k_{t+1}$, $x=\tilde{w}^k$, and $z=w^k_{t}$, 
	\begin{alignat}{6} 
	&&(1+\beta)||w^k_{t+1}-w^k_{t}||^2_2&\geq  
	||w^k_{t+1}-\tilde{w}^k||^2_2-\left(1+\frac{1}{\beta}\right)||w^k_{t}-\tilde{w}^k||^2_2\nonumber\\
	&&-\frac{13L_{E\lambda}}{4}||w^k_{t+1}-w^k_{t}||^2_2&\leq  
	-\frac{13L_{E\lambda}}{4(1+\beta)}||w^k_{t+1}-\tilde{w}^k||^2_2+\frac{13L_{E\lambda}\left(1+\frac{1}{\beta}\right)}{4(1+\beta)}||w^k_{t}-\tilde{w}^k||^2_2.\nonumber
	\end{alignat} 	
	Setting $\beta=2t-1$, 	
	$$-\frac{13L_{E\lambda}}{4}||w^k_{t+1}-w^k_{t}||^2_2\leq 
	-\frac{13L_{E\lambda}}{8t}||w^k_{t+1}-\tilde{w}^k||^2_2+\frac{13L_{E\lambda}}{8t-4}||w^k_{t}-\tilde{w}^k||^2_2.$$
	Applying this bound in \eqref{byi2},	
	$$\EE\h_{\lambda}(w^k_{t+1})\leq 	
	\EE\left[\h_{\lambda}(w^k_t)+\left(\frac{L_{E\lambda}}{6b}+\frac{13L_{E\lambda}}{8t-4}\right)||\tilde{w}^k-w^k_t||^2_2-
	\frac{13L_{E\lambda}}{8t}||w^k_{t+1}-\tilde{w}^k||^2_2-\frac{1}{36L_{E\lambda}}||\nabla 
	E^k_{t\lambda}(w^k_t)||^2_2\right].$$
	Summing over $t$,
	\begin{alignat}{6} 
	&\EE\h_{\lambda}(w^k_{m+1})&\leq&  
	\EE\left[\h_{\lambda}(w^k_1)+\sum_{t=1}^m\left(\frac{L_{E\lambda}}{6b}+\frac{13L_{E\lambda}}{8t-4}\right)||\tilde{w}^k-w^k_t||^2_2\right.\nonumber\\
	&&&\left.-\sum_{t=1}^{m}\frac{13L_{E\lambda}}{8t}||w^k_{t+1}-\tilde{w}^k||^2_2-\frac{1}{36L_{E\lambda}}\sum_{t=1}^m||\nabla
	E^k_{t\lambda}(w^k_t)||^2_2\right].\nonumber
	\end{alignat}	
	Considering that $\tilde{w}^k=w^k_1$ and $||w^k_{m+1}-\tilde{w}^k||^2_2\geq 0$,
	\begin{alignat}{6} 			
	&\EE\h_{\lambda}(w^k_{m+1})&\leq&\EE\left[  
	\h_{\lambda}(w^k_1)+\sum_{t=2}^m\left(\frac{L_{E\lambda}}{6b}+\frac{13L_{E\lambda}}{8t-4}\right)||\tilde{w}^k-w^k_t||^2_2\right.\nonumber\\
	&&&\left.-\sum_{t=1}^{m-1}\frac{13L_{E\lambda}}{8t}||w^k_{t+1}-\tilde{w}^k||^2_2-\frac{1}{36L_{E\lambda}}\sum_{t=1}^m||\nabla
	E^k_{t\lambda}(w^k_t)||^2_2\right]\nonumber\\	
	&&=&  
	\EE\left[\h_{\lambda}(w^k_1)+\sum_{t=1}^{m-1}\left(\frac{L_{E\lambda}}{6b}+\frac{13L_{E\lambda}}{8t+4}-\frac{13L_{E\lambda}}{8t}\right)||w^k_{t+1}-\tilde{w}^k||^2_2
	-\frac{1}{36L_{E\lambda}}\sum_{t=1}^m||\nabla 
	E^k_{t\lambda}(w^k_t)||^2_2\right]\nonumber\\	
	&&\leq &  
	\EE\left[\h_{\lambda}(w^k_1)+\sum_{t=1}^{m-1}\left(\frac{L_{E\lambda}}{6b}-\frac{L_{E\lambda}}{2t^2}\right)||w^k_{t+1}-\tilde{w}^k||^2_2-\frac{1}{36L_{E\lambda}}\sum_{t=1}^m||\nabla
	E^k_{t\lambda}(w^k_t)||^2_2\right]\nonumber\\	
	&&\leq & 
	\EE\left[\h_{\lambda}(w^k_1)-\frac{1}{36L_{E\lambda}}\sum_{t=1}^m||\nabla 
	E^k_{t\lambda}(w^k_t)||^2_2\right],\nonumber
	\end{alignat}		
	where the last inequality holds since $6b=6m^2>2(m-1)^2\geq 2t^2$ for $t=1,...,m-1$. This 
	summation can be equivalently written as 	
	\begin{alignat}{6} 
	&&\EE\h_{\lambda}(\tilde{w}^{k+1})&\leq 
	\EE\h_{\lambda}(\tilde{w}^k)-\EE\left[\frac{1}{36L_{E\lambda}}\sum_{t=1}^m||\nabla 
	E^k_{t\lambda}(w^k_t)||^2_2\right]\nonumber\\
	&&\EE\left[\frac{1}{36L_{E\lambda}}\sum_{t=1}^m||\nabla 
	E^k_{t\lambda}(w^k_t)||^2_2\right]&\leq 
	\EE\h_{\lambda}(\tilde{w}^k)-\EE\h_{\lambda}(\tilde{w}^{k+1})\nonumber\\
	&&\EE\left[\frac{1}{36L_{E\lambda}}\sum_{k=1}^S\sum_{t=1}^m||\nabla 
	E^k_{t\lambda}(w^k_t)||^2_2\right]&\leq 
	\h_{\lambda}(\tilde{w}^1)-\EE\h_{\lambda}(\tilde{w}^{S+1})\nonumber\\
	&&&\leq \h_{\lambda}(\tilde{w}^1)-\h_{\lambda}(w^*_{\lambda})\nonumber\\
	&&\EE\left[||\nabla 
	E^R_{T\lambda}(w^R_T)||^2_2\right]&\leq\frac{36L_{E\lambda}\left(\h_{\lambda}(\tilde{w}^1)-\h_{\lambda}(w^*_{\lambda})\right)}{Sm}.\nonumber\\
	&&&=\tilde{\DD}\frac{L+(Sm)^{\theta}}{Sm}.\nonumber
	\end{alignat}					
\end{proof}	

\section{Implementation details of SSD-SPG and SSD-SVRG}

In this section we describe all chosen parameter values using the notation found in 
\citep{xu2018}. The algorithm SSDC-SPG calls a stochastic proximal gradient (SPG) algorithm K times. 
For the $k^{th}$ iteration, the number of iterations of SPG equals $T_k=4k$. Each iteration of 
SPG uses one gradient call. We used the minimum $K$ which ensured at least $en$ gradient calls were 
used. The convex majorant parameter $\gamma=3L$, and the step size $\eta_t=1/(L(t+1))$. The Moreau 
envelope parameter $\mu=\epsilon$, where $K=O(1/\epsilon^4)$, is the only non-explicitly given 
parameter, which we set to $\mu=1/\left(K^{\frac{1}{4}}\right)$. SSDC-SVRG calls a stochastic variance 
reduced gradient (SVRG) algorithm $K$ times. We set the inner loop length $T_k=\max(2,200L/\gamma)$, 
and the outer loop length $S_k=\lceil \log_2(k)\rceil$. The step size $\eta_k=0.05/L$. Two 
parameters are not explicitly given, similar to in SSDC-SPG, we set 
$\mu=1/\left(K^{\frac{1}{4}}\right)$. For these parameter settings, there 
seems to be no restriction on $\gamma$. Their SVRG algorithm is based off of the work of 
\citet{xiao2014}, where empirical testing of different sizes of $T_k$ was done for a binary 
classification problem. The best performance was found with a choice of $T_k=2n$, from which we were 
able to determine $\gamma$. Given $\gamma$, we were then able to solve for $K$, ensuring at least $en$ 
gradient calls were used.

\bibliography{SSGM}
\bibliographystyle{icml2019}